\documentclass[12pt]{article}
\usepackage{amsfonts}
\usepackage{amssymb}
\usepackage{mathrsfs}
\usepackage{amsmath}
\usepackage{amsthm}
\usepackage{cite}

  \date{ }
     \title{{\bf Generalized  Derivations of Lie triple systems}
   \thanks{Supported by  NNSF of China (No. 11171055 and No. 11471090),  NSF of  Jilin province (No. 201115006).  Corresponding author (L. Chen):
     chenly640@nenu.edu.cn}}
   \author{Jia Zhou$^1$, Liangyun Chen$^2$, Yao Ma$^2$
     \vspace{0.3cm} \\$^1$School of Information Technology, Jilin Agricultural University,\\ Changchun, 130118, CHINA \\$^2$School of Mathematics and Statistics, Northeast Normal University,\\Changchun, 130024, CHINA}
 
\begin{document}
 \maketitle
 \begin{center}{\bf Abstract}\end{center}


In this paper, we present some basic properties concerning the
derivation algebra ${\rm Der}(T)$, the quasiderivation algebra ${\rm
QDer}(T)$ and the generalized derivation algebra ${\rm GDer}(T)$ of
a Lie triple system $T$, with the relationship ${\rm
Der}(T)\subseteq {\rm QDer}(T)\subseteq {\rm GDer}(T)\subseteq {\rm
End}(T)$. Furthermore, we completely determine those Lie triple
systems $T$ with condition ${\rm QDer}(T)={\rm End}(T)$. We also
show that the quasiderivations of $T$ can be embedded as derivations
in a larger Lie triple system.
   \vspace{0.3cm}

    \vspace{0.3cm}


   \vspace{0.3cm}
  \noindent {\bf \S 0\quad Introduction}
 \vspace{0.5cm}

  Lie triple systems arose initially in Cartan's study
   of Riemannian geometry. Jacobson {\rm\cite{J}} first introduced them in connection
with problems from Jordan theory and quantum mechanics, viewing Lie
triple systems as subspaces of Lie algebras that are closed related
to the ternary product. Lister gave the structure theory of
Lie triple systems of characteristic $0$ in {\rm\cite{L4}}. Hopkin introduced the concepts of
nilpotent ideals and the nil-radical of Lie triple systems, she
successfully generalized Engel's theorem to Lie triple systems in
characteristic zero{\rm\cite{H1}}. More recently, Lie triple systems
have been connected with the study of the Yang-Baxter equations
{\rm\cite{K}}.

As is well known, derivation and generalized derivation algebras are
very important subjects in the research of Lie algebras. In the
study of Levi factors in derivation algebras of nilpotent Lie
algebras, the generalized derivations, quasiderivations, centroids
and quasicentroids play key roles {\rm\cite{B}}. In {\rm\cite{M}},
Melville dealt particularly with the centroids of nilpotent Lie
algebras. The most important and systematic research on the
generalized derivation algebras of a Lie algebra and their
subalgebras was due to Leger and Luks. In {\rm\cite{L1}}, some nice
properties of  the quasiderivation algebras and of the
centroids have been obtained. In particular, they investigated the
structure of the generalized derivation algebras and characterized
the Lie algebras satisfying certain conditions. Meanwhile, they also
pointed that there exist some connections between quasiderivations
and cohomology of Lie algebras. For the generalized derivation
algebras of more general nonassociative algebras, the readers will
be referred to the papers\cite{CMN,F,F1,F2,H,JI,M,V}.

In this paper, we generalize some beautiful results in
{\rm\cite{L1}} to Lie triple system. In particular, we seek to
understand the structure of the generalized derivation algebras of a
Lie triple system or conversely, we want to characterize the Lie
triple systems for which the generalized derivation algebras or
their Lie subalgebras satisfy some special conditions.

This paper is organized as follows. Section 2 contains some
elementary observations about generalized derivations,
quasiderivations, centroids and quasicentroids, some of which are
technical results to be used in the sequel. In Section 3, we
characterize completely those Lie triple systems $T$ for which ${\rm
QDer}(T)={\rm End}(T)$. Such Lie triple systems include
two-dimensional simple Lie triple systems and all the commutative
Lie triple systems. Section 4 is devoted to showing that the
quasiderivations of a Lie triple system can be embedded as
derivations in a larger Lie triple system ${\breve{T}}$. Moreover,
if the center of $T$ is zero, we obtain a semidirect decomposition
of ${\rm Der}(\breve{T})$. \vspace{0.3cm}

  \noindent {\bf \S 1\quad Preliminaries}
 \vspace{0.5cm}

 \noindent{\bf Definition 1.1} {\rm\cite{M1}} {\it A Lie triple system is a pair $(T,[\cdot,\cdot,\cdot])$ consisting of a vector space $T$
over a field $\mathbb{F}$, a trilinear multiplication
$[\cdot,\cdot,\cdot]:T\times T\times T\rightarrow T$ such that for
all $x, y, z, u, v \in T,$
$$\begin{array}{ll}[x, x, z]=0,\end{array}$$
$$\begin{array}{ll}[x, y, z]+ [y,z,x]+[z,x,y]=0,\end{array}$$
$$\begin{array}{ll}[x,y,[z,u,v]]=[[x,y,z],u,v]+[z,[x,y,u],v]+[z,u,[x,y,v]].\end{array}$$

\vspace{0.3cm}

{\it ${\rm End}(T)$ denotes the set consists of all linear maps of
$T$. Obviously, ${\rm End}(T)$ is a Lie algebra over $\mathbb{F}$
with the bracket $[D_{1},D_{2}]=D_{1}D_{2}-D_{2}D_{1},$ for all
$D_{1},D_{2}\in {\rm End}(T).$}\vspace{0.3cm}

\noindent{\bf Definition 1.2} {\rm\cite{S}} Let {\it
$(T,[\cdot,\cdot,\cdot])$ be a Lie triple system. A linear map
$D:T\rightarrow T$ is said to be a derivation of $T$ if it satisfies
$$[D(x),y,z]+[x,D(y),z]+[x,y,D(z)]=D([x,y,z]),$$
$\forall x,y,z\in T.$}

We denote the set of all derivations by ${\rm Der}(T)$, then ${\rm
Der}(T)$ provided with the commutator is a subalgebra of ${\rm
End}(T)$ and is called the derivation algebra of $T$.\vspace{0.3cm}

 \noindent{\bf Definition 1.3} {\it $D\in {\rm End}(T)$ is said to be a generalized
 derivation of $T$, if there exist
$D',D'',D'''\in {\rm End}(T)$ such that
$$[D(x),y,z]+[x,D^{'}(y),z]+[x,y,D^{''}(z)]=D^{'''}([x,y,z]),\eqno(1.1)$$
for all $x, y,z\in T.$} \vspace{0.3cm}

 \noindent{\bf Definition 1.4} {\it $D\in {\rm
End}(T)$ is said to be a quasiderivation, if there exists $D'\in
{\rm End}(T)$ such that
$$[D(x),y,z]+[x,D(y),z]+[x,y,D(z)]=D^{'}([x,y,z]),\eqno(1.2)$$ for all $x, y,z\in
T.$}

Denote by ${\rm GDer}(T)$ and ${\rm QDer}(T)$ the sets of
generalized derivations and quasiderivations,
respectively.\vspace{0.3cm}

 \noindent{\bf Definition 1.5} {\rm\cite{L3}} {\it If ${\rm
 C}(T)=\{D\in {\rm End}(T)\mid[D(x),y,z]=[x,D(y),z]=[x,y,D(z)]\\=D([x,y,z])\}$
for all $x,y,z\in T,$ then ${\rm C}(T)$ is called a centroid of
$T$.}\vspace{0.3cm}

 \noindent{\bf Definition 1.6} If {\it
${\rm QC}(T)=\{D\in {\rm End}(T)\mid[D(x),y,z]=[x,D(y),z]=[x,y,D(z)]
\}$ for all $x,y,z\in T, $ then ${\rm QC}(T)$ is called a
quasicentroid of $T$.}\vspace{0.3cm}

\noindent{\bf Definition 1.7} If {\it ${\rm ZDer}(T)=\{D\in {\rm
End}(T)\mid[D(x),y,z]=D([x,y,z])=0 \}$ for all $x,y,z\in T, $ then
${\rm ZDer}(T)$ is called a central derivation of $T$.}

It is easy to verify that $${\rm ZDer}(T)\subseteq {\rm
Der}(T)\subseteq {\rm QDer}(T)\subseteq {\rm GDer}(T)\subseteq {\rm
End}(T).$$

\noindent{\bf Definition 1.8} {\rm\cite{L3}} {\it $T$ is a Lie
triple system and $I$ is a non-empty subset of $T$. We call ${\rm
Z}_{T}(I)=\{x\in T| [x,a,y]=[y,a,x]=0, \forall a\in I, y\in T\}$ the
centralizer of $I$ in $T$. In particular, ${\rm Z}_{T}(T)=\{x\in T|
[x,y,z]=0, \forall y,z\in T\}$ is the center of $T$, denoted by
${\rm Z}(T)$}.\vspace{0.3cm}

     \noindent{\bf \S 2\quad Generalized derivation algebras and their  subalgebras}
 \vspace{0.5cm}

  First, we  give some basic properties
of center derivation algebra, quasiderivation algebra and the
generalized derivation algebra of a Lie triple system.

     \noindent{\bf  Proposition 2.1} {\it Let
$T$ be a Lie triple system. Then the following statements hold:

$(1)$\quad ${\rm GDer}(T),{\rm QGer}(T)$ and ${\rm C}(T)$ are
subalgebras of  ${\rm End}(T)$.

$(2)$\quad ${\rm ZDer}(T)$ is an ideal of ${\rm Der}(T)$.}

\rm {\it Proof.}\quad  $(1)$ Assume that $D_{1},D_{2}\in{\rm
GDer}(T)$. For all $x,y,z\in T,$ we have
\begin{eqnarray*}
[D_{1}D_{2}(x),y,z]&=&D^{'''}_{1}D^{'''}_{2}([x,y,z])-D^{'''}_{1}[x,D^{'}_{2}(y),z]-D^{'''}_{1}[x,y,D^{''}_{2}(z)]\\
&-&[D_{2}(x),D^{'}_{1}(y),z]-[D_{2}(x),y,D^{''}_{1}(z)]\\
&=&D^{'''}_{1}D^{'''}_{2}([x,y,z])-[D_{1}(x),D^{'}_{2}(y),z]-[x,D^{'}_{1}D^{'}_{2}(y),z] \\
&-&[x,D^{'}_{2}(y),D^{''}_{1}(z)]-[D_{1}(x),y,D^{''}_{2}(z)]-[x,D^{'}_{1}(y),D^{''}_{2}(z)]\\
&-&[x,y,D^{''}_{1}D^{''}_{2}(z)]-[D_{2}(x),D^{'}_{1}(y),z]-[D_{2}(x),y,D^{''}_{1}(z)],
\end{eqnarray*}and

\begin{eqnarray*}
[D_{2}D_{1}(x),y,z]&=&D^{'''}_{2}D^{'''}_{1}([x,y,z])-D^{'''}_{2}[x,D^{'}_{1}(y),z]-D^{'''}_{2}[x,y,D^{''}_{1}(z)]\\
&-&[D_{1}(x),D^{'}_{2}(y),z]-[D_{1}(x),y,D^{''}_{2}(z)]\\
&=&D^{'''}_{2}D^{'''}_{1}([x,y,z])-[D_{2}(x),D^{'}_{1}(y),z]-[x,D^{'}_{2}D^{'}_{1}(y),z]\\
&-&[x,D^{'}_{1}(y),D^{''}_{2}(z)]-[D_{2}(x),y,D^{''}_{1}(z)]-[x,D^{'}_{2}(y),D^{''}_{1}(z)]\\
&-&[x,y,D^{''}_{2}D^{''}_{1}(z)]-[D_{1}(x),D^{'}_{2}(y),z]-[D_{1}(x),y,D^{''}_{2}(z)]
\end{eqnarray*}

\noindent Thus for all $x,y,z\in T$, we
have$$[[D_{1},D_{2}](x),y,z]=[D'''_{1},D'''_{2}]([x,y,z])-[x,y,[D''_{1},D''_{2}](z)]-[x,[D'_{1},D'_{2}](y),z].$$

From the definition of generalized derivation, one gets
$[D_{1},D_{2}]\in {\rm GDer}(T),$ so ${\rm GDer}(T)$ is a subalgebra
of ${\rm End}(T)$.

Similarly, ${\rm QGer}(T)$ is a subalgebra of ${\rm End}(T)$.

Assume that $D_{1},D_{2}\in{\rm C}(T).~\forall x,y,z\in T,$ note
that
\begin{eqnarray*}
[[D_{1},D_{2}](x),y,z]&=&[D_{1}D_{2}(x),y,z]-[D_{2}D_{1}(x),y,z]\\
&=&D_{1}([D_{2}(x),y,z])-D_{2}([D_{1}(x),y,z])\\
&=&D_{1}D_{2}([x,y,z])-D_{2}D_{1}([x,y,z])\\
&=&[D_{1},D_{2}]([x,y,z]).
\end{eqnarray*}
Similarly,
$$[x,[D_{1},D_{2}](y),z]=[D_{1},D_{2}]([x,y,z])=[x,y,[D_{1},D_{2}](z)].$$
Then $[D_{1},D_{2}]\in {\rm C}(T)$, ${\rm C}(T)$ is a subalgebra of
${\rm End}(T)$.

$(2)$ Assume that $D_{1}\in{\rm ZDer}(T),D_{2}\in{\rm Der}(T).$ For
all $x,y,z\in T,$ we have
$$\begin{array}{ll}[[D_{1},D_{2}]([x,y,z])]=D_{1}D_{2}([x,y,z])-D_{2}D_{1}([x,y,z])=0,\end{array}$$
and
\begin{eqnarray*}
[[D_{1},D_{2}](x),y,z]&=&[(D_{1}D_{2}-D_{2}D_{1})(x),y,z]\\
&=&D_{1}([D_{2}(x),y,z])-[D_{1}(x),D_{2}(y),z]=0.
\end{eqnarray*}
Then $[D_{1},D_{2}]\in {\rm ZDer}(T)$ and ${\rm ZDer}(T)$ is an
ideal of ${\rm Der}(T)$.\hfill$\Box$\vspace{0.3cm}

\noindent{\bf Lemma 2.2} {\it Let $T$ be a Lie triple system. Then

    $(1)$ \quad $[{\rm Der}(T),{\rm C}(T)]\subseteq {\rm C}(T);$

    $(2)$ \quad $[{\rm QDer}(T),{\rm QC}(T)]\subseteq {\rm QC}(T);$

    $(3)$ \quad ${\rm C}(T)\cdot{\rm Der}(T)\subseteq {\rm Der}(T);$

    $(4)$ \quad ${\rm C}(T)\subseteq {\rm QDer}(T);$

    $(5)$ \quad $[{\rm QC}(T),{\rm QC}(T)]\subseteq {\rm QDer}(T);$

    $(6)$ \quad ${\rm QDer}(T)+{\rm QC}(T)\subseteq {\rm GDer}(T).$}

 \rm {\it Proof.}\quad $(1)$ Assume that
$D_{1}\in{\rm Der}(T),D_{2}\in{\rm C}(T).$ For all $x,y,z\in T,$ we
have
\begin{eqnarray*}
[D_{1}D_{2}(x),y,z]&=&D_{1}([D_{2}(x),y,z])-[D_{2}(x),D_{1}(y),z]-[D_{2}(x),y,D_{1}(z)]\\
&=&D_{1}D_{2}([x,y,z])-[x,D_{2}D_{1}(y),z]-[x,y,D_{2}D_{1}(z)],
\end{eqnarray*}and
\begin{eqnarray*}
[D_{2}D_{1}(x),y,z]&=&D_{2}(D_{1}([x,y,z])-[x,D_{1}(y),z]-[x,y,D_{1}(z)])\\
&=&D_{2}D_{1}([x,y,z])-[x,D_{2}D_{1}(y),z]-[x,y,D_{2}D_{1}(z)].
\end{eqnarray*}Hence,
\begin{align*}[[D_{1},D_{2}](x),y,z]&=D_{1}D_{2}([x,y,z])-D_{2}D_{1}([x,y,z])=[D_{1},D_{2}]([x,y,z]).\end{align*}
 Similarly,
\begin{align*}[[D_{1},D_{2}](x),y,z]=[x,[D_{1},D_{2}](y),z]=[x,y,[D_{1},D_{2}](z)].\end{align*}
Thus, $[D_{1},D_{2}]\in{\rm C}(T)$ and we get $[{\rm Der}(T),{\rm
C}(T)]\subseteq {\rm C}(T).$

$(2)$ Similar to the proof of $(1)$.

$(3)$ Assume that $D_{1}\in{\rm C}(T),D_{2}\in{\rm Der}(T).$ For all
$x,y,z\in T,$ we have
\begin{eqnarray*}
D_{1}D_{2}[x,y,z]&=&D_{1}([D_{2}(x),y,z])+[x,D_{2}(y),z]+[x,y,D_{2}(z)])\\
&=&[D_{1}D_{2}(x),y,z]+[x,D_{1}D_{2}(y),z]+[x,y,D_{1}D_{2}(z)].
\end{eqnarray*}So we have $D_{1}D_{2}\in {\rm Der}(T).$

$(4)$ Assume that $D\in{\rm QC}(T).$ For all $x,y,z\in T,$ we have
$$\begin{array}{ll}[D(x),y,z]=[x,D(y),z]=[x,y,D(z)].\end{array}$$
Hence,
$$\begin{array}{ll}[D(x),y,z]+[x,D(y),z]+[x,y,D(z)]=3D[x,y,z].\end{array}$$
Therefore, $D\in {\rm QDer}(T)$ since $D^{'}=3D\in {\rm
C}(T)\subseteq {\rm End}(T)$.

$(5)$ Assume that $D_{1},D_{2}\in{\rm QC}(T).$ For all $x,y,z\in T,$
we have\\

$~~~~~[[D_{1},D_{2}](x),y,z]+[x,[D_{1},D_{2}](y),z]+[x,y,[D_{1},D_{2}](z)]
\\~~~~~~~~=[D_{1}D_{2}(x),y,z]+[x,D_{1}D_{2}(y),z]+[x,y,D_{1}D_{2}(z)]-[D_{2}D_{1}(x),y,z]
\\~~~~~~~~-[x,D_{2}D_{1}(y),z]-[x,y,D_{2}D_{1}(z)].$

\noindent And
$$\begin{array}{ll}[D_{1}D_{2}(x),y,z]&=[D_{2}(x),D_{1}(y),z]=[x,D_{2}D_{1}(y),z],\end{array}$$
$$\begin{array}{ll}[D_{1}D_{2}(x),y,z]&=[D_{2}(x),y,D_{1}(z)]=[x,y,D_{2}D_{1}(z)].\end{array}$$
\noindent
Hence,$$\begin{array}{ll}[[D_{1},D_{2}](x),y,z]+[x,[D_{1},D_{2}](y),z]+[x,y,[D_{1},D_{2}](z)]=0,\end{array}$$
i.e. $[D_{1},D_{2}]\in {\rm QDer}(T)$.

$(6)$\quad It is obvious.\hfill$\Box$
 \vspace{0.3cm}

\noindent{\bf Lemma 2.3} {\rm\cite{S}} {\it If $T$ is a Lie triple
system, $I$ is an ideal of $T$, then ${\rm Z}_{T}(I)$ is also an
ideal of $T$. Moreover, ${\rm Z}(T)={\rm Z}_{T}(T), {\rm Z}(I)={\rm
Z}_{I}(I)$ are ideals of $T$.
 \vspace{0.3cm}

\noindent{\bf Lemma 2.4} {\rm\cite{S}} {Let the Lie triple system
$T$ be decomposed into the direct sum of two ideals, i.e. $T=A\oplus
B$. Then we have

$(1)\quad {\rm Z}(T)={\rm Z}(A)\oplus{\rm Z}(B)$.

$(2)$\quad If ${\rm Z}(T)=0,$ then ${\rm Der}(T)={\rm
Der}(A)\oplus{\rm Der}(B)$.
 \vspace{0.3cm}

\noindent{\bf Proposition 2.5} {If the Lie triple system $T$ can be
decomposed into the direct sum of two ideals, i.e. $T=A\oplus B$ and
${\rm Z}(T)=0,$ then we have

$(1)$\quad ${\rm GDer}(T)={\rm GDer}(A)\oplus{\rm GDer}(B)$;

$(2)$\quad ${\rm QDer}(T)={\rm QDer}(A)\oplus{\rm QDer}(B)$;

$(3)$\quad ${\rm C}(T)={\rm C}(A)\oplus{\rm C}(B)$;

$(4)$\quad ${\rm QC}(T)={\rm QC}(A)\oplus{\rm QC}(B)$.

\rm {\it Proof.}\quad $(1)$\quad For $D^{'}\in {\rm GDer}(A),$
extend it to a linear transformation on $T$ by setting
$D^{'}(a+b)=D^{'}(a), \forall a\in A,b\in B.$ Obviously, $D^{'}\in
{\rm GDer}(T)$ and ${\rm GDer}(A)\subseteq {\rm GDer}(T)$.
Similarly, ${\rm GDer}(B)\subseteq {\rm GDer}(T).$ Let $a\in A,
b_{1},b_{2}\in B$ and $D\in {\rm Der}(T).$
Then $$\begin{array}{ll}[D(a),b_{1},b_{2}]&=D([a,b_{1},b_{2}])-[a,D(b_{1}),b_{2}]-[a,b_{1},D(b_{2})]\\
&=-[a,D(b_{1}),b_{2}]-[a,b_{1},D(b_{2})]\in A\cap
B=\{0\}.\end{array}$$ Suppose $D(a)=a^{'}+b^{'},$ where $a\in A,
b\in B,$
then$$\begin{array}{ll}0=[D(a),b_{1},b_{2}]&=[a^{'},b_{1},b_{2}]+[b^{'},b_{1},b_{2}].\end{array}$$
So $[b^{'},b_{1},b_{2}]=0$ and $b^{'}\in Z(B).$ Since ${\rm
Z}(T)={\rm Z}(A)\oplus{\rm Z}(B), b^{'}=0.$ Hence $D(a)=a^{'}\in A.$
Therefore $D(A)\subseteq A.$ Similarly, $D(B)\subseteq B.$

Let $D\in {\rm Der}(T)$ and $x=a+b\in A+B,$ where $a\in A,~b\in B$.
Define $E, F\in {\rm End}(T)$ by $E(a+b)=D(a), F(a+b)=D(b),$ then
$E\in {\rm Der}(A), F\in {\rm Der}(B)$. Hence $D=E+F\in {\rm
Der}(A)+{\rm Der}(B).$ Since ${\rm Der}(A)\cap{\rm Der}(B)=\{0\},
{\rm GDer}(T)={\rm GDer}(A)\dotplus{\rm GDer}(B)$ as a vector space.

Let $E\in {\rm Der}(A), F\in {\rm Der}(B)$ and $b\in B.$ Then
$[E,D]=(ED-DE)(b)=0.$ Hence $[E,D]\in {\rm Der}(A)$ and ${\rm
Der}(A)\triangleleft {\rm Der}(T).$ Similarly, ${\rm
Der}(B)\triangleleft {\rm Der}(T).$

$(2),(3),(4)$\quad Similar to the proof of $(1)$. \hfill$\Box$
 \vspace{0.3cm}

\noindent{\bf Proposition 2.6} {\it If $T$ is a Lie triple system,
then ${\rm QC}(T)+[{\rm QC}(T),{\rm QC}(T)]$ is a subalgebra of
${\rm GDer}(T)$.}

{\it Proof.}\quad By the conclutions of Lemma $2.2~(5)$ and $(6)$,
we have
$${\rm QC}(T)+[{\rm QC}(T),{\rm QC}(T)]\subseteq {\rm GDer}(T)$$ and
\begin{eqnarray*}
&&\ \ [{\rm QC}(T)+[{\rm QC}(T),{\rm QC}(T)],{\rm QC}(T)+[{\rm QC}(T),{\rm QC}(T)]]\\
&&\subseteq[{\rm QC}(T)+{\rm QDer}(T),{\rm QC}(T)+[{\rm QC}(T),{\rm QC}(T)]]\\
&&\subseteq[{\rm QC}(T),{\rm QC}(T)]+[{\rm QC}(T),[{\rm QC}(T),{\rm
QC}(T)]]+[{\rm QDer}(T),{\rm QC}(T)]\\&&+[{\rm QDer}(T),[{\rm
QC}(T),{\rm QC}(T)]].\end{eqnarray*} It is easy to verify $[{\rm
QDer}(L),[{\rm QC}(L),{\rm QC}(L)]]\subseteq [{\rm QC}(L),{\rm
QC}(L)]$ by the Jacobi identity of  Lie algebra. Thus $$[{\rm
QC}(L)+[{\rm QC}(L),{\rm QC}(L)],{\rm QC}(L)+[{\rm QC}(L),{\rm
QC}(L)]]\subseteq {\rm QC}(L)+[{\rm QC}(L),{\rm QC}(L)].$$
 \vspace{0.3cm}

\noindent{\bf Theorem 2.7} {\it If $T$ is a Lie triple system, then
$[{\rm C}(T),{\rm QC}(T)]\subseteq {\rm End}(T,{\rm Z}(T)).$
 Moreover, if ${\rm
Z}(T)=\{0\},$ then $[{\rm C}(T),{\rm QC}(T)]=\{0\}.$}

{\it Proof.}\quad Assume that $D_{1}\in {\rm C}(T), D_{2}\in {\rm
QC}(T)$ and for all $x,y,z\in T$, we have
$$\begin{array}{ll}[[D_{1},D_{2}](x),y,z]&=[D_{1}D_{2}(x),y,z]-[D_{2}D_{1}(x),y,z]\\
&=D_{1}([D_{2}(x),y,z])-[D_{1}(x),D_{2}(y),z]\\
&=D_{1}([D_{2}(x),y,z]-[x,D_{2}(y),z])=0.\end{array}$$ Hence
$[D_{1},D_{2}](x)\in {\rm Z}(T)$ and $[D_{1},D_{2}]\in {\rm
End}(T,{\rm Z}(T))$ as desired. Furthermore, if ${\rm Z}(T)=\{0\},$
it is clear that $[{\rm C}(T),{\rm QC}(T)]=\{0\}.$ \hfill$\Box$
 \vspace{0.3cm}

\noindent{\bf  Definition 2.8} {\rm\cite{Z}} {\it Let $L$ be an
algebra over $\mathbb F$ (char $\mathbb F\not=2$), if the
multiplication satisfies the following identities:
$$x\cdot y=y\cdot x, $$
$$(((x\cdot y)\cdot w)\cdot z-(x\cdot y)\cdot(w\cdot z))+(((y\cdot z)\cdot w)\cdot x-
(y\cdot z)\cdot(w\cdot x))$$
$$+(((z\cdot x)\cdot w)\cdot y-(z\cdot x)\cdot(w\cdot y))=0,$$
for all $x,y,z,w\in L,$ then we call $L$ a Jordan algebra.}
\vspace{0.3cm}

\noindent{\bf  Proposition 2.9} {\rm\cite{Z}} \it Let $T$ be a Lie
triple system over $\mathbb F$ (char $\mathbb F\not=2$), with the
operation $D_{1}\bullet D_{2}=D_{1}D_{2}+D_{2}D_{1},~$ for all
elements $D_{1},D_{2}\in {\rm End}(T)$.
 Then the pair $({\rm End}(T),\bullet)$ is a Jordan
algebra.} \vspace{0.3cm}

\noindent{\bf  Corollary 2.10} \it Let $T$ be a Lie triple system
over $\mathbb F$ (char $\mathbb F\not=2$), with the operation
$D_{1}\bullet D_{2}=D_{1}D_{2}+D_{2}D_{1},~$ for all elements
$D_{1},D_{2}\in {\rm QC}(T)$.
 Then $({\rm QC}(T)$ is a Jordan
algebra.}

\rm {\it Proof.}\quad We need only to show that $D_{1}\bullet
D_{2}\in{\rm QC}(T)$, for all $D_{1},D_{2}\in{\rm QC}(T),$ we have
\begin{eqnarray*}
[D_{1}\bullet D_{2}(x),y,z]&=&[D_{1}D_{2}(x),y,z]+[D_{2}D_{1}(x),y,z]\\
&=&[D_{2}(x),D_{1}(y),z]+[D_{1}(x),D_{2}(y),z]\\
&=&[x,D_{2}D_{1}(y),z]+[x,D_{1}D_{2}(y),z]\\
&=&[x,D_{1}\bullet D_{2}(y),z].
\end{eqnarray*}
Similarly, $[D_{1}\bullet D_{2}(x),y,z]=[x,y,D_{1}\bullet
D_{2}(z)]$. Then $D_{1}\bullet D_{2}\in {\rm QC}(T)$ and ${\rm
QC}(T)$ is a Jordan algebra.

\vspace{0.3cm}

\noindent{\bf Theorem 2.11} {\it If $T$ is a Lie triple system over
$\mathbb F$, then we have

$(1)$\quad If char $\mathbb F\not=2$, then ${\rm QC}(T)$ is a Lie
algebra with $[D_{1},D_{2}]= D_{1}D_{2}-D_{2}D_{1}$ if and only if
${\rm QC}(T)$ is also an associative algebra with respect to
composition.

$(2)$\quad If char $\mathbb F\not=3$ and ${\rm Z}(T)=\{0\}$, then
${\rm QC}(T)$ is a Lie algebra if and only if $~[{\rm QC}(T),{\rm
QC}(T)]=0.$}

\rm {\it Proof.}\quad $(1)$ $(\Leftarrow)$ For all $D_{1},D_{2}\in
{\rm QC}(T),$ we have $D_{1}D_{2}\in {\rm QC}(T)$ and $D_{2}D_{1}\in
{\rm QC}(T)$, so $[D_{1},D_{2}]=D_{1}D_{2}-D_{2}D_{1}\in {\rm
QC}(T).$ Hence, ${\rm QC}(T)$ is a Lie algebra.

$(\Rightarrow)$ Note that $D_{1}D_{2}=D_{1}\bullet
D_{2}+\frac{[D_{1},D_{2}]}{2}$ and by Corollary $2.10$,  we have
$D_{1}\bullet D_{2}\in {\rm QC}(T), [D_{1},D_{2}]\in {\rm QC}(T)$.
It follows that $D_{1}D_{2}\in {\rm QC}(T)$ as desired.

$(2)$ $(\Rightarrow)$\quad Assume that $D_{1},~D_{2}\in {\rm
QC}(T).$ For all $x,y,z\in T$, ${\rm QC}(T)$ is a Lie algebra, so
$[D_{1},D_{2}]\in {\rm QC}(T)$, then
$$[[D_{1},D_{2}](x),y,z]=[x,[D_{1},D_{2}](y),z]=[x,y,[D_{1},D_{2}](z)].$$
From the proof of Lemma$~2.2~(5)$, we have
$$[[D_{1},D_{2}](x),y,z]=-[x,[D_{1},D_{2}](y),z]-[x,y,[D_{1},D_{2}](z)].$$
Hence $3[[D_{1},D_{2}](x),y,z]=0.$ We have
$[[D_{1},D_{2}](x),y,z]=0,$ i.e. $[D_{1},D_{2}]=0$ since char
$\mathbb F\not=3$.

$(\Leftarrow)$\quad It is clear.\hfill$\Box$
    \vspace{0.3cm}

\noindent{\bf Lemma 2.12} {\it Let $V$ be a linear space and
$\mathcal {A}:V \rightarrow V$ a linear map. $f(x)$ denotes the
minimal polynomial of $f$. If $x^{2}$ does not divide $f(x)$, then
$V={\rm Ker}(\mathcal {A})\dotplus {\rm Im}(\mathcal {A}).$

{\it Proof.}\quad \rm Obviously ${\rm dim}(V)={\rm dim}({\rm
Ker}(\mathcal {A}))+{\rm dim}({\rm Im}(\mathcal {A}))$ because
$\mathcal {A}$ is a linear map. $x^{2}$ does not divide $f(x)$ means
that $f(x)=x^{2}g(x)+ax+b,~a\not=0$ or $b\not=0.$

Case $1$: If $b\not=0$, then $f(\mathcal {A})=\mathcal
{A}^{2}g(\mathcal {A})+a\mathcal {A}+b\rm id =0,$ so $\mathcal
{A}(\mathcal {A}g(\mathcal {A})+a\rm id)={\it-b}id$ and $\mathcal
{A}$ is invertible. Hence ${\rm Ker}(\mathcal {A})=\{0\}.$

Case $2$: If $b=0,$ that means $a\not=0,~f(\mathcal {A})=\mathcal
{A}^{2}g(\mathcal {A})+a\mathcal {A}.$ Here we only prove ${\rm
Ker}(\mathcal {A})\cap {\rm Im}(\mathcal {A})=\{0\}$. Indeed,
$\forall x\in {\rm Ker}(\mathcal {A})\cap {\rm Im}(\mathcal
{A}),~\mathcal {A}(x)=0$ and there exists $x^{'}\in V$ such that
$\mathcal {A}(x^{'})=x.$ So $f(\mathcal {A})(x^{'})=\mathcal
{A}^{2}g(\mathcal {A})(x^{'})+a\mathcal {A}(x^{'}),$ which means
$a\mathcal {A}(x^{'})=ax=0$. Hence $x=0$ since
$a\not=0.$\hfill$\Box$
 \vspace{0.3cm}

\noindent{\bf Proposition 2.13} {\it Let $T$ be a Lie triple system,
$D\in {\rm C}(T)$. Then

$(1)$\quad ${\rm Ker}(D)$and ${\rm Im}(D)$ are ideals in $T$.

$(2)$\quad If $T$ is indecomposable, $D\in {\rm C}(T)$ and
$D\not=0.$ Suppose $x^{2}$ does not divide the minimal polynomial of
$D$, then $D$ is invertible.

$(3)$\quad If $T$ is indecomposable and ${\rm C}(T)$ consists of
semisimple elements, then ${\rm C}(T)$ is a field.

\rm {\it Proof.} $(1)$\quad Since $D\in {\rm C}(T),$ for all $x\in
{\rm Ker}(D),~y,z\in T,$ one gets $D[x,y,z]=[D(x),y,z]=0,$ that
means $[x,y,z]\in {\rm Ker}(D)$.

Meanwhile, for all $x\in {\rm Im}(D),$ there exists $x'\in T,$ such
that $x=D(x')$. So $[x,y,z]=[D(x'),y,z]=D[x',y,z]\in {\rm Im}(D)$.

$(2)$\quad From Lemma {\rm 2.12} and $(1)$ there is an ideal sum
$T={\rm Ker}(D)\oplus {\rm Im}(D).$ Since $T$ is indecomposable, one
gets ${\rm Ker}(D)=0$ and ${\rm Im}(D)=T,$ which means $D$ is
invertible.

$(3)$\quad For all semisimple element $D\in {\rm C}(T),$ since
$x^{2}$ does not divide the minimal polynomial of $D$ and $T$ is
indecomposable, from $(2)$ one gets $D$ is invertible. It is obvious
that ${\rm id}\in {\rm C}(T).$ If there exist
$D_{1}\not=0,~D_{2}\not=0,~D_{1},D_{2}\in {\rm C}(T)$ such that
$D_{1}D_{2}=0,$ then $D_{1}=D_{2}=0$, a contradiction. Hence ${\rm
C}(T)$ has no zero divisor. Obviously, one gets
$D_{1}D_{2}=D_{2}D_{1},~\forall D_{1},D_{2}\in {\rm C}(T).$ So ${\rm
C}(T)$ is a field. \hfill$\Box$ \vspace{0.3cm}

\noindent{\bf  Lemma 2.14} {\it Let $T$ be a Lie triple system with
${\rm Z}(T)=\{0\}$. If $D\in {\rm QC}(T)$ and suppose $x^{2}$ does
not divide the minimal polynomial of $D$, then $T={\rm Ker}(D)\oplus
{\rm Im}(D).$}

{\it Proof.}\quad From Lemma {\rm 2.12}, there is a vector space
direct sum $T={\rm Ker}(D)\dotplus {\rm Im}(D).$ Obviously, $[{\rm
Ker}(D),D(T),T]=[D({\rm Ker}(D)),T,T]=0$ and $[T,D(T),{\rm
Ker}(D)]\\=[T,T,D({\rm Ker}(D))]=0,$ so ${\rm Ker}(D)\subseteq {\rm
Z}_{T}({\rm Im}(D)),~{\rm Im}(D)\subseteq {\rm Z}_{T}({\rm
Ker}(D)).$ Since ${\rm Z}_{T}({\rm Im}(D))\cap {\rm Z}_{T}({\rm
Ker}(D))={\rm Z}(T)=\{0\},$ we must have ${\rm Ker}(D)={\rm
Z}_{T}({\rm Im}(D)),~{\rm Im}(D)\\={\rm Z}_{T}({\rm Ker}(D)).$

It is easy to get $[[{\rm Ker}(D),T,T],{\rm Im}(D),T]=[T,{\rm
Im}(D),[{\rm Ker}(D),T,T]]=0$, which means $[{\rm
Ker}(D),T,T]\subseteq {\rm Z}_{T}({\rm Im}(D))={\rm Ker}(D).$

Also $[[{\rm Im}(D),T,T],{\rm Ker}(D),T]=[T,{\rm Ker}(D),[{\rm
Im}(D),T,T]]=0,~[{\rm Im}(D),T,T]\subseteq {\rm Z}_{T}({\rm
Ker}(D))={\rm Im}(D).$ So ${\rm Ker}(D)$ and ${\rm Im}(D)$ are
ideals. \hfill$\Box$ \vspace{0.3cm}

\noindent{\bf  Corollary 2.15} {\it Let $(T,[\cdot,\cdot,\cdot])$ be
a indecomposable Lie triple system over an algebraically field
$\mathbb F$ and $~{\rm Z}(T)=0$. $D\in {\rm QC}(T)$ is semisimple,
then $D\in Z_{{\rm C}(T)}({\rm GDer}(T)).$}

{\it Proof.}\quad Let $D\in {\rm QC}(T),$ $D$ has an eigenvalue
$\lambda$ since $\mathbb F$ is an algebraically field. We denote the
corresponding eigenspace by $E_{\lambda}(D),$ it is easy to get
$(D-\lambda {\rm id})\in {\rm QC}(T)$ and ${\rm Ker}(D-\lambda {\rm
id})=E_{\lambda}(D)\not=0.$ From Lemma $2.14$ and $D$ is a
semisimple element, ${\rm Ker}(D-\lambda {\rm id})$ is an ideal of
$T$, so ${\rm Ker}(D-\lambda {\rm id})=T.$ That is $D=\lambda {\rm
id}\in {\rm C}(T)$ and $[D,{\rm GDer}(T)]=0.$ \hfill$\Box$
\vspace{0.3cm}

\noindent {\bf \S 3\quad  Lie triple systems with ${\rm
QDer}(T)={\rm End}(T)$}
        \vspace{0.3cm}

Let $T$ be a Lie triple system over $\mathbb{F}$. We define a linear
map $\phi$ : $T\otimes T\otimes T\rightarrow\ T, x\otimes y\otimes
z\mapsto[x,y,z]$. Define ${\rm Ker}(\phi)$:=$\{ \sum x\otimes
y\otimes z\in T\otimes T\otimes T\mid x,y,z \in {T},~\sum
[x,y,z]=0\}$, then it is easy to see that ${\rm Ker}(\phi)$
 is a subspace of $T\otimes T\otimes
T$.

We define $(T\otimes\ T\otimes\ T)^{+}:=\langle x\otimes y\otimes
z+y\otimes x\otimes z\mid x,y,z \in T\rangle$ and $(T\otimes\
T\otimes\ T)^{-}:=\langle x\otimes y\otimes z-y\otimes x\otimes
z\mid x,y,z \in \ T\rangle$, then both $(T\otimes\ T\otimes\ T)^{+}$
and $(T\otimes\ T\otimes\ T)^{-}$ are subspaces of $T\otimes
T\otimes T$. It is easy to check that
$$T\otimes T\otimes T=(T\otimes\ T\otimes\ T)^{+}\dotplus(T\otimes\ T\otimes\ T)^{-} (\textrm{direct
sum of vector spaces}),$$ and we also have $\textrm{dim}(T\otimes\
T\otimes\ T)^{+}=n^{2}(n+1)/2$ and $\textrm{dim}(T\otimes\ T\otimes\
T)^{-}=n^{2}(n-1)/2,$ where $\textrm{dim}(T)=n.$

For all $D\in {\rm End}(T)$, we define $D^{\ast}\in {\rm
End}(T\otimes T\otimes T)$ satisfying that
$$D^{\ast}(x\otimes y\otimes z)=D(x)\otimes
y\otimes z+x\otimes D(y)\otimes z+x\otimes y\otimes
D(z),\eqno(3.1)$$ for all $x,y,z\in T$.

\noindent{\bf Lemma 3.1} $D\in {\rm QDer}(T)$ \textit{if and only
if} $D^{\ast}({\rm Ker}(\phi))\subseteq {\rm Ker}(\phi)$.

\textit{Proof.}\quad $(\Rightarrow)$ For all $\sum x\otimes y\otimes
z\in {\rm Ker}(\phi)$, we have $\sum [x,y,z]=0$.
Thus,\begin{eqnarray*}
D^{\ast}(\sum x\otimes y\otimes z)&=&\sum D^{\ast}(x\otimes y\otimes z)\\
&=&\sum (D(x)\otimes y\otimes z+x\otimes D(y)\otimes z+x\otimes
y\otimes D(z)).
\end{eqnarray*}Since $D\in {\rm QDer}(T)$, we have

$$\sum ([D(x),y,z]+[x,
D(y),z]+[x,y,D(z)])=\sum D'([x,y,z])=D'(\sum[x,y,z])=0.$$ Hence
$D^{\ast}({\rm Ker}(\phi))\subseteq {\rm Ker}(\phi)$.

$(\Leftarrow)$ Since $D^{\ast}({\rm Ker}(\phi))\subseteq {\rm
Ker}(\phi)$, there exists an element $D'\in {\rm End}(T)$ such that
$$\phi\circ D^{\ast}=D'\circ \phi:\ T\otimes T\otimes T\rightarrow\ T,$$ for
all $D\in {\rm End}(T)$. A direct computation shows that for all
$x,y,z\in T$,
$$[D(x),y,z]+[x,
D(y),z]+[x,y,D(z)]=D'([x,y,z]),$$ that is, $D\in {\rm
QDer}(T)$.\hfill$\Box$\vspace{0.3cm}

\noindent{\bf Lemma 3.2} \textit{Suppose that} ${\rm
End}(T)$\textit{ acts on} $T\otimes T\otimes T$ \textit{via}
$D\cdot(x\otimes y\otimes z)=D^{\ast}(x\otimes y\otimes z)$
\textit{for all} $x,y,z\in T$ \textit{with} $D^{\ast}$ \textit{as
above. Then }$(T\otimes\ T\otimes\ T)^{+}$ and $(T\otimes\ T\otimes\
T)^{-}$ \textit{are two irreducible} ${\rm
End}(T)$-\textit{modules.}

\textit{Proof.}\quad Here we need only to show that $(T\otimes\
T\otimes\ T)^{+}$ is an irreducible ${\rm End}(T)$-module and the
case of $(T\otimes\ T\otimes\ T)^{-}$ is similar. For all $x,y,z\in
T$ and $D\in {\rm End}(T)$, from Eq.$(3.1)$ we have
\begin{eqnarray*}
&&\ \ D\cdot\left(x\otimes y\otimes z+y\otimes x\otimes z\right)\\
&&=D\cdot(x\otimes y\otimes z)+D\cdot(y\otimes x\otimes z)\\
&&=D(x)\otimes y\otimes z+x\otimes D(y)\otimes z+x\otimes y\otimes
D(z)+D(y)\otimes x\otimes z\\
&&\quad +y\otimes D(x)\otimes z+y\otimes x\otimes
D(z)\\
&&=(D(x)\otimes y\otimes z+y\otimes D(x)\otimes z)+(x\otimes y\otimes D(z)+y\otimes x\otimes D(z))\\
&&\quad +(x\otimes D(y)\otimes z+D(y)\otimes x\otimes z)\in
(T\otimes\ T\otimes\ T)^{+}.\end{eqnarray*} Hence $D\cdot(T\otimes\
T\otimes\ T)^{+} \subseteq(T\otimes\ T\otimes\ T)^{+}$ and it is
easy to check that
$$[D_{1},D_{2}]\cdot(x\otimes y\otimes z)=D_{1}\cdot(
D_{2}\cdot(x\otimes y\otimes z))-D_{2}\cdot( D_{1}\cdot (x\otimes
y\otimes z)).$$ Therefore, $(T\otimes\ T\otimes\ T)^{+}$ is an ${\rm
End}(T)$-module.

Suppose that there exists a nonzero ${\rm End}(T)$-submodule $V$ in
$(T\otimes\ T\otimes\ T)^{+}$. Choose a nonzero element $\sum
(x\otimes y\otimes z+y\otimes x\otimes z)$ in $V$ for some $x,y,z
\in T$.  A direct computation shows that all elements in $(T\otimes\
T\otimes\ T)^{+}$ are obtainable by repeated application of elements
of ${\rm End}(T)$ to $\sum (x\otimes y\otimes z+y\otimes x\otimes
z)$ and formation of linear combinations.  Hence $V$ is $(T\otimes\
T\otimes\ T)^{+}$ itself. Thus, $(T\otimes\ T\otimes\ T)^{+}$ as an
${\rm End}(T)$-module is irreducible.
\hfill$\Box$\\

\noindent{\bf Theorem 3.3} \textit{Let $T$ be a Lie triple system
with $[T,T,T]\not=0$ and} ${\rm QDer}(T)={\rm End}(T)$.
\textit{Then} $T$ \textit{is a two-dimensional simple Lie triple
system. }

\textit{Proof.}\quad We consider the action of ${\rm End}(T)$ on
$T\otimes\ T\otimes\ T$ via $D\cdot (x\otimes y\otimes
z)=D^{\ast}(x\otimes y\otimes z)$ for all $x,y,z\in T$ with
$D^{\ast}$ as in Eq.$~(3.1)$. By Lemma $3.1$, ${\rm QDer}(T)={\rm
End}(T)$ implies that ${\rm End}(T)\cdot {\rm Ker}(\phi)\subseteq
{\rm Ker}(\phi)$. Lemma $3.2$ tells us that the only proper
subspaces of $T\otimes\ T\otimes\ T$, invariant under this action of
${\rm End}(T)$, are $(T\otimes\ T\otimes\ T)^{+}$ and $(T\otimes\
T\otimes\ T)^{-}$. Thus we have $\phi:T\otimes\ T\otimes\
T\rightarrow\ T$ with kernel $\{0\},(T\otimes\ T\otimes\ T)^{+}$ and
$(T\otimes\ T\otimes\ T)^{-}$. Using
$${\rm dim}(T)\geq {\rm dim}(T\otimes\ T\otimes\ T)-{\rm dim(Ker}(\phi)),$$ we have that $n=1$, if
${\rm Ker}(\phi)=\{0\}$ or $(T\otimes\ T\otimes\ T)^{-}$; $n\leq2$,
if ${\rm Ker}(\phi)=(T\otimes\ T\otimes\ T)^{+}$.

We discuss the possibilities for $n$ as follows:

(\textit{a}) If $n=1$, then $T$ is commutative, which is a
contradiction with our assumption.

(\textit{b}) If $n=2$, i.e.$~{\rm Ker}(\phi)=(T\otimes\ T\otimes\
T)^{+}$, hence ${\rm dim(Ker}(\phi))=6$ and ${\rm dim}([T,T,T])=2$,
so $\phi$ must be surjective and  we have $[T,T,T]=T$. So that $T$
is the two-dimensional simple Lie triple system (See
{\rm\cite[Chapter $4.3$]{L2}}). \hfill$\Box$ \vspace{0.3cm}

 On the other hand, the
converse of Theorem $3.3$ is also valid. One can prove the following
theorem.

\noindent{\bf Theorem 3.5} \textit{If $T$ is a two-dimensional
simple Lie triple system or an abelian Lie triple system, then}
${\rm QDer}(T)={\rm End}(T)$.

\textit{Proof.}\quad By {\rm\cite[Chapter $4.3$]{L2}}, we have a
basis $e_{1}, e_{2}$ such that $[e_{1},e_{2},e_{1}]=-e_{1},$
$[e_{1},e_{2},e_{2}]=e_{2}.$ Then for all $D\in {\rm End}(T),$ one
gets $D(e_{1})=k_{1}e_{1}+k^{'}_{2}e_{2},~k_{1},~k^{'}_{1},$
$k_{1},~k^{'}_{2}\in \mathbb F.$ It obvious $k^{'}=0$ since $D$ is a
linear map. Similarly, $D(e_{2})=k_{2}e_{2}$. So we have
$$[D(e_{1}),e_{2},e_{1}]+[e_{1},D(e_{2}),e_{1}]+[e_{1},e_{2},D(e_{1})]=-(2k_{1}+k_{2})e_{1}.$$
$$[D(e_{1}),e_{2},e_{2}]+[e_{1},D(e_{2}),e_{2}]+[e_{1},e_{2},D(e_{2})]=(k_{1}+2k_{2})e_{2}.$$
Let $D^{'}\in {\rm End}(T)$ such that
$D^{'}(e_{1})=-(2k_{1}+k_{2})e_{1}$ and
$D^{'}(e_{2})=(k_{1}+2k_{2})e_{2},$ thus $D\in {\rm QDer}(T)$.
  \hfill$\Box$

\vspace{0.3cm}

\noindent {\bf \S 4\quad  The quasiderivations of Lie triple
systems}
        \vspace{0.3cm}

       In this section, we will prove that the quasiderivations of $T$ can
be embedded as derivations in a larger Lie triple system and obtain
a direct sum decomposition of {\rm Der}($T$) when the center ${\rm
Z}(T)$ is equal to zero.

 \noindent{\bf Proposition 4.1} {\it Let $T$ be a Lie triple system over
       ${\mathbb F}$ and $t$ an indeterminate. We define $\breve{T}:=
\{\Sigma(x\otimes t+y\otimes t^{3})| x,y\in T\}$. Then $\breve{T}$
is a Lie triple system with the operation $[x\otimes t^{i},y\otimes
t^{j},z\otimes t^{k}]=[x,y,z]\otimes t^{i+j+k},$ for all $x,y,z\in
T,i,j,k\in\{1,3\}$}.

{\it Proof.}\quad For all $x,y,z,u,v\in T$ and $i,j,k,m,n
\in\{1,3\},$ we have
$$\begin{array}{ll}[x\otimes t^{i},y\otimes t^{j},z\otimes t^{k}]&=[x,y,z]\otimes t^{i+j+k}\\
&=-[y,x,z]\otimes t^{i+j+k}\\
&=-[y\otimes t^{j},x\otimes t^{i},z\otimes t^{k}],\end{array}$$
\begin{eqnarray*}
&&\ \ [x\otimes t^{i},y\otimes t^{j},z\otimes t^{k}]+[y\otimes t^{j},z\otimes t^{k},x\otimes t^{i}]+[z\otimes t^{k},x\otimes t^{i},y\otimes t^{j}]\\
&&=[x,y,z]\otimes t^{i+j+k}+[y,z,x]\otimes t^{i+j+k}+[z,x,y]\otimes t^{i+j+k}\\
&&=([x,y,z]+[y,z,x]+[z,x,y])\otimes t^{i+j+k}=0,\end{eqnarray*}
and\begin{eqnarray*}
&&\ \ [x\otimes t^{i},y\otimes t^{j},[z\otimes t^{k},u\otimes t^{m},v\otimes t^{n}]]=[x,y,[z,u,v]]\otimes t^{i+j+k+m+n}\\
&&=([[x,y,z],u,v]+[z,[x,y,u],v]+[z,u,[x,y,v]])\otimes t^{i+j+k+m+n}\\
&&=[[x\otimes t^{i},y\otimes t^{j},z\otimes t^{k}],u\otimes t^{m},v\otimes t^{n}]+[z\otimes t^{k},[x\otimes t^{i},y\otimes t^{j},u\otimes t^{m}],v\otimes t^{n}]\\
&&+[z\otimes t^{k},u\otimes t^{m},[x\otimes t^{i},y\otimes
t^{j},v\otimes t^{n}]].\end{eqnarray*} Hence $\breve{T}$ is a Lie
triple system. \hfill$\Box$ \vspace{0.3cm}

For convenience, we write $xt(xt^{3})$ in place of $x\otimes
t(x\otimes t^{3}).$

If $U$ is a subspace of $T$ such that $T=U\oplus [T,T,T],$ then
$$\breve{T}=Tt+Tt^{3}=Tt+Ut^{3}+[T,T,T]t^{3},$$
Now we define a map $\varphi:{\rm QDer}(T)\rightarrow {\rm
End}(\breve{T})$ satisfying
$$\varphi(D)(at+ut^{2}+bt^{2})=D(a)t+D'(b)t^{3},$$
where $D\in {\rm QDer}($T$),$ and $D'$ is in Eq.$(1.2)$, $a\in
T,u\in U,b\in [T,T,T]$. \vspace{0.3cm}

  \noindent{\bf Proposition 4.2} {\it $T,\breve{T},\varphi$ are as defined above.Then

 $(1)$ $\varphi$ is injective and $\varphi(D)$ does not depend
on the choice of $D'$.

$(2)$ $\varphi({\rm QDer}(T))\subseteq {\rm Der}(\breve{T}).$}

{\it Proof.}\quad (1)\quad If $\varphi(D_{1})=\varphi(D_{2}),$ then
for all $a\in T,b\in [T,T,T]$ and $u\in U,$ we have
$$\varphi(D_{1})(at+ut^{3}+bt^{3})=\varphi(D_{2})(at+ut^{3}+bt^{3}),$$ that is $$D_{1}(a)t+D'_{2}(b)t^{3}=
D_{2}(a)t+D'_{2}(b)t^{3},$$ so $D_{1}(a)=D_{2}(a).$ Hence
$D_{1}=D_{2},$ and $\varphi$ is injective.

Suppose that there exists $D''$ such that
$$\varphi(D)(at+ut^{3}+bt^{3})=D(a)t+D''(b)t^{3},$$ and
$$[D(x),y,z]+[x,D(y),z]+[x,y,D(z)]=D''([x,y,z]),$$
then we have $$D'([x,y,z])=D''([x,y,z]),$$ thus $D'(b)=D''(b).$
Hence
$$\varphi(D)(at+ut^{3}+bt^{3})=D(a)t+D'(b)t^{3}=D(a)t+D''(b)t^{3},$$
which implies $\varphi(D)$ is determined by $D$.

(2)\quad We have $[xt^{i},yt^{j},zt^{k}]=[x,y,z]t^{i+j+k}=0,$ for
all $i+j+k\geq 4$. Thus, to show $\varphi(D)\in {\rm
Der}(\breve{T}),$ we need only to check the validness of the
following equation
$$\varphi(D)([xt,yt,zt])=[\varphi(D)(xt),yt,zt]+[xt,\varphi(D)(yt),zt]+[xt,yt,\varphi(D)(zt)].$$
For all $x,y,z\in T,$ we have
$$\begin{array}{ll}\varphi(D)([xt,yt,zt])&=\varphi(D)([x,y,z]t^3)=D'([x,y,z])t^3\\
&=([D(x),y,z]+[x,D(y),z]+[x,y,D(z)])t^3\\
&=[D(x)t,yt,zt]+[xt,D(y)t,zt]+[xt,yt,D(z)t]\\
&=[\varphi(D)(xt),yt,zt]+[xt,\varphi(D)(yt),zt]+[xt,yt,\varphi(D)(zt)].\end{array}$$
Therefore, for all $D\in {\rm QDer}(T)$, we have $\varphi(D)\in {\rm
Der}(\breve{T})$.\hfill$\Box$ \vspace{0.3cm}

\noindent{\bf Proposition 4.3} {\it Let $T$ be a Lie triple system.
${\rm Z}(T)=\{0\}$ and $\breve{T},~\varphi$ are as defined above.
Then ${\rm Der}(\breve{T})=\varphi({\rm QDer}(T))\dotplus {\rm
ZDer}(\breve{T}).$}

{\it Proof.}\quad Since ${\rm Z}(T)=\{0\}$, we have ${\rm
Z}(\breve{T})=Tt^3.$ For all $g\in {\rm Der}(\breve{T}),$ we have
$g({\rm Z}(\breve{T}))\subseteq {\rm Z}(\breve{T}),$ hence
$g(Ut^3)\subseteq g({\rm Z}(\breve{T}))\subseteq {\rm
Z}(\breve{T})=Tt^3.$ Now we define a map
$f:Tt+Ut^3+[T,T,T]t^3\rightarrow Tt^3$ by
$$\ f(x)=\left\{\begin{array}{ll}g(x)\cap Tt^3,& x\in Tt ;\\
 g(x),& x\in Ut^3 ;\\  0,& x\in [T,T]t^3.\end{array}\right.$$
It is clear that $f$ is linear. Note that
$$f([\breve{T},\breve{T},\breve{T}])=f([T,T,T]t^3)=0,$$
$$~[f(\breve{T}),\breve{T},\breve{T}]\subseteq [Tt^3,Tt+Tt^3,Tt+Tt^3]=0,$$ hence $f\in {\rm ZDer}(\breve{T}).$
Since $$(g-f)(Tt)=g(Tt)-g(Tt)\cap Tt^{3}=g(Tt)-Tt^{3}\subseteq Tt,~
(g-f)(Ut^3)=0,$$ and
$$(g-f)([T,T,T]t^3)=g([\breve{T},\breve{T},\breve{T}])\subseteq
[\breve{T},\breve{T},\breve{T}]=[T,T,T]t^3,$$ there exist $D,~D'\in
{\rm End}(T)$ such that for all $a\in T,~b\in [T,T,T]$,
$$(g-f)(at)=D(a)t,~ (g-f)(bt^3)=D'(b)t^3.$$ Since $(g-f)\in {\rm
Der}(\breve{T})$ and by the definition of ${\rm Der}(\breve{T})$, we
have
$$[(g-f)(a_1t),a_2t,a_3t]+[a_1t,(g-f)(a_2t),a_3t]+[a_1t,a_2t,(g-f)(a_3t)]=(g-f)([a_1t,a_2t,a_3t]),$$
for all $a_1,a_2,a_3\in T.$ Hence
$$[D(a_1),a_2,a_3]+[a_1,D(a_2),a_3]+[a_1,a_2,D(a_3)]=D'([a_1,a_2,,a_3]).$$ Thus
$D\in {\rm QDer}(T).$ Therefore, $g-f=\varphi(D)\in \varphi({\rm
QDer}(T))$, so ${\rm Der}(\breve{T})\subseteq \varphi({\rm
QDer}(T))+{\rm ZDer}(\breve{T}).$ By Proposition $4.2~(2)$ we have
${\rm Der}(\breve{T})=\varphi({\rm QDer}(T))+{\rm ZDer}(\breve{T}).$

For all $f\in \varphi({\rm QDer}(T))\cap{\rm ZDer}(\breve{T})$,
there exists an element $D\in {\rm QDer}(T)$ such that
$f=\varphi(D).$ Then
$$f(at+ut^3+bt^3)=\varphi(D)(at+ut^3+bt^3)=D(a)t+D'(b)t^3,$$ for all $a\in T,b\in [T,T,T].$

On the other hand, since $f\in {\rm ZDer}(\breve{T}),$ we have
$$f(at+bt^3+ut^3)\in {\rm Z}(\breve{T})=Tt^3.$$ That is to say,
$D(a)=0,$ for all $a\in T$ and so $D=0.$ Hence $f=0.$

Therefore ${\rm Der}(\breve{T})=\varphi({\rm QDer}(T))\dotplus {\rm
ZDer}(\breve{T})$ as desired. \hfill$\Box$


\vspace{0.3cm}

\begin{thebibliography}{99}

\bibitem{B} Benoist, Y. (1988). La partie semi-simple de l'alg$\grave{\textrm{e}}$bre des d$\acute{\textrm{e}}$rivations
d'une alg$\grave{\textrm{e}}$bre de Lie nilpotente. C.R. Acad. Sci.
Paris. 307:901-904.
\bibitem{CMN} Chen L. Y., Ma Y., Ni L. (2013). Generalized Derivations of Lie Color Algebras. Results Math. 63:923-936.
\bibitem{F} Fialkow, L. A. (1980). Generalied derivations. Topics in Modern Operator Theory. 95-103. Operator Theory: Adv. Appl., 2, Birkh\"{a}user, Basel-Boston, Mass., 1981.
\bibitem{F1} Filippov, V. T. (2000). On $\delta$-derivations of prime alternative and Mal'tsev algebras. (Russian)Algebra Log.
39:618-625.
\bibitem{F2} Filippov, V. T. (1998). On $\delta$-derivations of Lie algebras. (Russian). Sibirsk. Mat. Zh. 39:1409-1422.
\bibitem{H} Hopkins, N. C. (1996). Generalized derivations of nonassociative algebras. Nova J. Math. Game Theory Algebra. 5:215-224.
\bibitem{H1} Hopkins. N. C. (1985). Some structure theory for a class of triple systems.
Trans, Amer. Math. Soc. 289:203-212.
\bibitem{J} Jacobson N. (1949). Lie and Jordan triple Systems. Amer.
J. Math. Soc. 71:149-170.
\bibitem{JI} Jimenez-Gestal, C., Perez-Izquierdo, J. M. (2008). Ternary derivations of finite-dimensional real division
algebras. Linear Algebra Appl. 428:2192-2219.
\bibitem{K} Kamiya N., Okubo S. (1997). On triple systems and
Yang-Baxter equations//Proceedings of the Seventh International
Colloquium on Differential Equations. Utrecht:VSP. 189-196.
\bibitem{L1} Leger, G. F., Luks, E. M. (2000).  Generalized derivations of Lie algebras. J. Algebra. 228:165-203.
\bibitem{L2} Lev S., Larissa S. and Ivan S. (2006). Non-Associative Algebra and Its Applications. CRC Press, U.S.A.
\bibitem{L3} Lin, J. (2010). Centroid of Lie triple systems. Acta Scientiarum Naturalium Universitatis Nankaiensis. 43:98-104.
\bibitem{L4} Lister, William G.  (1952). A structure theory of Lie Triple
systems. Trans, Amer. Math. Soc. 72:217-242.
\bibitem{M} Melville, D. J. (1992). Centroids of nilpotent Lie algebras. Comm. Algebra. 20:3649-3682.
\bibitem{M1} Ma, Y.,Chen, L. Y., Lin, J. (2014). Systems of quotients of Lie triple systems. Comm. Algebra. 42:3339-3349.
\bibitem{S} Shi, Y. Q., Meng, D. J. (2002). On derivations and automorphism group of Lie triple systems. Acta Scientiarum Naturalium Universitatis Nankaiensis. 35:32-37.
\bibitem{V} Vinberg, \`{E}. B. (1989). Generalized derivations of
algebras. Algebra and analysis. 185-188. Amer. Math. Soc. Transl. Ser. 2, 163, Amer. Math. Soc., Providence, RI, 1995.
\bibitem{Z} Zhang, R. X., Zhang, Y. Z. (2010). Generalized derivations of Lie superalgebras. Comm. Algebra. 38:3737-3751.

\end{thebibliography}
\end{document}